\documentclass[12pt]{koval}

\usepackage{eucal}

\def\D{\widehat D(A)}
\def\wtilde{\widetilde}
\def\tPsi{\wtilde\Psi}
\def\DC{\widehat D(\tPsi)}
\def\DB{{\mathcal D}(\tPsi)}
\def\SC{{\mathcal D}'(\RE^n)}
\def\SCT{\widetilde{\mathcal D}'(\RE^n)}
\def\E{{\mathcal E}'(\RE^n)}
\def\B{\mathsf B}
\def\L{\mathsf L}
\def\H{\mathcal H}
\def\X{\mathcal X}
\def\Y{\mathcal Y}
\def\G{\mathcal G}
\def\F{\mathcal F}
\def\RE{\mathbb R}
\def\C{{\mathbb C}}
\def\LD{L^2(\RE^n)}
\def\phireg{\phi_{\text{\rm reg}}}
\def\vp{\varphi}
\def\vpreg{\vp_{\text{\rm reg}}}
\def\-{\,\text{\rm -}\,}
\def\h{\text{\rm h}}
\def\p{\par\noindent}

\begin{document}
\p
Operator Theory:\p
Advances and Application, Vol. 132, 333-346
\vskip 20pt
\title[Boundary conditions for Singular Perturbations]{Boundary
Conditions for Singular Perturbations of Self-Adjoint Operators}

\author{Andrea Posilicano}

\address{Dipartimento di Scienze, Universit\`a dell'Insubria, I-22100
Como, Italy}

\email{andreap@uninsubria.it}

\begin{abstract}
 Let $A:D(A)\subseteq\H\to\H$ be an injective
self-adjoint operator and let $\tau:D(A)\to\X$, $\X$ a Banach
space, be a surjective
linear map such that $\|\tau\phi\|_\X\le c\,\|A\phi\|_\H$. Supposing 
that \text{\rm Kernel}$\,\tau$ is dense in $\H$, we 
define a family $A^\tau_\Theta$ of self-adjoint operators which are
extensions of the symmetric operator $A_{\left|\left\{\tau=0\right\}\right.}$. 
Any $\phi$ in the operator domain $D(A^\tau_\Theta)$ is characterized by a sort
of boundary conditions on its univocally defined regular component
$\phireg$, which belongs to the completion of
$D(A)$ w.r.t. the norm $\|A\phi\|_\H$. These boundary conditions are written
in terms of the map $\tau$, playing the role of a trace (restriction)
operator, as $\tau\phireg=\Theta\, Q_\phi$, the extension parameter 
$\Theta$ being a self-adjoint operator
from $\X'$ to $\X$. The self-adjoint extension is then simply defined by
$A^\tau_\Theta\phi:=A\,\phireg$. 
The case in which $A\phi=\Psi*\phi$ is a convolution operator on
$\LD$, $\Psi$ a distribution with compact support, is studied in detail.
\end{abstract}

\maketitle

\section{Introduction}
Let $$A:D(A)\subseteq \H\to\H$$ be a self-adjoint operator on the complex Hilbert
space $\H$ (to prevent any misunderstanding we 
remark here that all over the paper we will avoid to identify a
Hilbert space with its strong dual). As usual $D(A)$
inherits a Hilbert space structure by introducing the scalar
product leading to the graph norm
$$
\|\phi\|_A^2:=\langle\phi,\phi\rangle_\H+\langle A\phi,A\phi\rangle_\H\,.
$$  
Considering then a linear bounded operator
$$
\tau:D(A)\to\X \,,\qquad \tau\in \B(D(A),\X)\ ,
$$
$\X$ a complex Banach space, we are interested in describing
the self-adjoint extensions of the symmetric operator $A_{\left|\{\tau=0\}\right.}$. 
In typical situations $A$ is a
(pseudo-)differential 
operator on $\LD$ and $\tau$ is a trace (restriction) operator along
some null 
subset $F\subset\RE^n$ (see e.g. \cite{[AFHKL]}-\cite{[AK]}, \cite{[Br]}-\cite{[C]},
\cite{[KKO]}-\cite {[P]}, \cite {[Sh]}, \cite{[T]} and
references therein). 
\par 
Denoting the 
resolvent set of $A$ by $\rho(A)$, we define $R(z)\in \B(\H,D(A))$,
$z\in\rho(A)$, by
$$
R(z):=(-A+z)^{-1}
$$
and we then introduce, for any $z\in\rho(A)$, the operators $\breve G(z)\in
\B(\H,\X)$ and $G(z)\in\wtilde \B(\X',\H)$ by
$$
\breve G(z):=\tau\cdot R(z)\,,\qquad 
G(z):=C_\H^{-1}\cdot\breve G(z^*)'
\ .\eqno(1)
$$
Here the prime $'$ denotes both the strong dual space and the
(Banach) adjoint map, 
and $C_\H$ indicates the canonical conjugate-linear isomorphism on $\H$ to
$\H'$ (the reader is refered to section 2 below for a list of
definitions and notations). As an immediate consequence of the first resolvent 
identity for $R(z)$ we have (see [{\bf 19}, lemma 2.1]) 
$$
(z-w)\, R(w)\cdot G(z)=G(w)-G(z)\eqno(2)
$$
and so 
$$
\forall\,w,z\in\rho(A),\qquad \text{\rm
Range}(G(w)-G(z))\subseteq D(A)\,.\eqno(3)
$$
In [{\bf 19}, thm. 2.1], by means of a Kre\u\i n-like formula, and under the hypotheses 
$$\tau\,\text{ is surjective}\eqno (\h 1) $$
$$\text{\rm Range}\,\tau'\cap\H'
=\left\{0\right\}\,,\eqno (\h2)$$
we constructed a
family $A^\tau_\Theta$ of self-adjoint extension of
$A_{\left|\left\{\tau=0\right\}
\right.}$ by giving its resolvent family. The hypothesis (h1) could be
weakened, see \cite{[P]}, but here we
prefer to use a simpler framework. In formulating (h2) we used the
embedding of $\H'$ into $D(A)'\supseteq\,$Range$\,\tau'$ given by 
$\varphi\mapsto \langle C_\H^{-1}\varphi,\,\cdot\,\rangle_\H$. 
Such an hypothesis is then equivalent to the denseness, in $\H$, of
the set $\left\{\tau=0\right\}$. Indeed there exists $\ell\in\X'$ such 
that $\tau' \ell\in\H'$ if and
only if there exists $\psi\in\H$ (necessarily orthogonal to Kernel$\,\tau$) such
that for any $\phi\in D(A)$ one has 
$\langle\psi,\phi\rangle_\H=\ell(\tau\phi)$.\par 
The advantage of the formula
given in \cite{[P]} over other approaches (see e.g. \cite{[S]}, 
\cite{[DM1]}, \cite{[DM2]}, \cite{[GMT]} and references therein) is its
relative simplicity, being expressed directly in terms of the map
$\tau$; moreover the domain of
definition of $A^\tau_\Theta$ can be described,
interpreting the map $\tau$ as a trace (restriction) operator, in terms
of a sort of boundary conditions (see [{\bf 19}, remark 2.10]). In the case 
$0\notin\sigma(A)$, $\sigma(A)$ denoting the spectrum of $A$, this description becomes particularly expressive 
since $A^\tau_\Theta\phi$ can be simply defined
by the original operator applied to the regular component of $\phi$. 
Such a regular component 
$\phi_0\in D(A)$ is univocally determined by the natural decomposition
which 
enter in the definition
of $D(A_\Theta^\tau)$ and it has to satisfy the boundary
condition 
$$\tau\phi_0=\Theta\,Q_\phi\,.$$ 
More precisely, by (h1), (h2) and by
[{\bf 19}, lemma 2.2, thm. 2.1, prop. 2.1, remarks 2.10, 2.12], we have the following
\begin{theorem} Let $A:D(A)\subseteq\H\to\H$ be self-adjoint with
$0\notin\sigma(A)$, let $\tau:D(A)\to\X$ be continuous and satisfy
{\rm (h1)} and {\rm (h2)}. If
$\Theta\in\wtilde \L(\X',\X)$ is self-adjoint, $G:=G(0)$ and
$$
D(A_\Theta^\tau):=\left\{\phi\in\H\, :\, \phi=
\phi_0+GQ_\phi,\, \phi_0\in
D(A),\, Q_\phi\in D(\Theta),\, \tau\phi_0=\Theta\,Q_\phi\right\},
$$
then the linear operator
$$
A_\Theta^\tau:D(A_\Theta^\tau)\subseteq\H\to\H\,,\qquad A_\Theta^\tau\phi:=A\phi_0
$$
is self-adjoint and
coincides with $A$ on the kernel of $\tau$; 
the decomposition entering in the definition of its domain is unique. 
Its resolvent is given by
$$
R^\tau_\Theta(z):=R(z)+G(z)\cdot(\Theta+\Gamma(z))^{-1}\cdot\breve
G(z)\,,
\qquad z\in W^-_\Theta\cup W_\Theta^+\cup\C\backslash\RE\, ,
$$
where $$\Gamma(z):=\tau\cdot(G-G(z))$$ and 
$$
W^\pm_\Theta:=\left\{\,\lambda\in\RE\cap \rho(A)\ :\ 
\gamma(\pm\Gamma(\lambda))>-\gamma(\pm\Theta)\,\right\}\ .
$$
\end{theorem}
\begin{remark} By (h1) one has $\X\simeq D(A)/\text{\rm Kernel$\,\tau$}
\simeq (\text{\rm Kernel$\,\tau$})^\perp$ and so 
$$D(A)\simeq \text{\rm Kernel$\,\tau$}\oplus\X\,.$$
This implies that $\X$ inherits a
Hilbert space structure and we could then identify $\X'$ with $\X$. 
Even if this gives some advantage (see
[{\bf 19}, remarks 2.13-2.16, lemma 2.4]) here we prefer to use only the
Banach space structure of $\X$.    
\end{remark}
The purpose of the present paper is to extend the above theorem to the
case in which $A$ is merely injective. 
Thus, denoting the pure point
spectrum of $A$ by $\sigma_{\text{\rm }pp}(A)$, we require $0\notin\sigma_{\text{\rm }pp}(A)$ but we do not
exclude the case $0\in\sigma(A)\backslash\sigma_{\text{\rm }pp}(A)$;
this is a typical situation when 
$A$ is a differential operator on $\LD$. In order to carry out this program 
we will suppose that the map $\tau$ has a continuos extension to
$\D$, the completion of $D(A)$ with respect to the norm
$\|A\phi\|_\H$ (note that  $\D=D(A)$ when $0\notin\sigma(A)$). 
This further hypothesis allows then to perform
the limit 
$\lim_{\epsilon\to 0}G(i\epsilon)-G(z)$ (see lemma 3); 
thus an analogue on
the above theorem 1 is obtained (see theorem 5). Such an abstract
construction is successively specialized to the case in
which $A\phi=\Psi*\phi$ is a convolution operator on $\LD$, where
$\Psi$ is a distribution with compact support (so that this comprises
the case of differential-difference operators). In this situation the results
obtained in theorem 5 can be made more appealing (see theorem 11). The
case in which $A=\Delta:H^2(\RE^n)\to \LD$, $n>4$, and $\tau$ is the trace
(restriction) operator along a $d$-set with a
compact closure of zero Lebesgue measure, $0<n-d<4$, is explicitly studied
(see example 14). Of course, since $-\Delta$ is not negative, in this
case one could apply theorem 1 to $-\Delta+\lambda$, $\lambda>0$, and then
define
$-\Delta^\tau_\Theta:=(-\Delta+\lambda)^\tau_\Theta-\lambda$. However 
this alternative definition looks a bit artificial and has
the drawback of giving rise to boundary conditions which depend on the
arbitrary parameter $\lambda$. 
The starting motivation of this work was indeed the desire to get rid
of such a dependence. 
\section{Definitions and notations}
\begin{itemize}
\item Given a Banach space $\X$ we denote by $\X'$ its strong dual;  
\item $\L(\X,\Y)$, resp. $\wtilde \L(\X,\Y)$, denotes the
space of linear, resp. conjugate linear, operators 
from the Banach space $\X$ to the Banach space $\Y$; 
$\L(\X):=\L(\X,\X)$, $\wtilde \L(\X):=\wtilde \L(\X,\X)$. 
\item $\B(\X,\Y)$, resp. $\wtilde \B(\X,\Y)$, denotes the (Banach) space of 
bounded, everywhere defined, linear, resp. conjugate linear, operators 
on the Banach space $\X$ to the Banach space $\Y$. 
\item Given $A\in \L(\X,\Y)$ and $\wtilde A\in\wtilde \L(\X,\Y)$
densely defined, the
closed operators $A'\in \L(\Y',\X')$ and $\wtilde A'\in\wtilde
\L(\Y',\X')$ the are the adjoints of $A$ and $\wtilde A$ respectively,
i.e. 
\begin{gather*}\forall\,x\in D(A)\subseteq\X,\quad\forall
\ell\in D(A')\subseteq
\Y',\qquad (A'\ell)(x)=\ell(Ax)\,,\\
\forall\,x\in D(\wtilde
A)\subseteq\X,\quad \forall\ell\in D(\wtilde
A')\subseteq \Y',\qquad(\wtilde A'\ell)(x)=(\,\ell(\wtilde Ax)\,)^*
\end{gather*} 
where $^*$ denotes complex conjugation. 
\item $J_\X\in \B(\X,\X'')$ indicates the injective map (an
isomorphism when $\X$ is reflexive) defined by $(J_\X\, x)(\ell):=\ell(x)$.
\item A closed, densely defined operator $A\in \L(\X',\X)\cup\wtilde \L(\X',\X)$
is said to be self-adjoint if $J_\X\cdot A=A'$. 
\item For any self-adjoint $A\in\L(\X',\X)\cup\wtilde \L(\X',\X)$ we
define 
$$\gamma(A):=
\inf\left\{\,\ell(A \ell),\ \ell\in
D(A),\ \|\ell\|_{\X'}=1\,\right\}\,.
$$
\item If $\H$ is a complex Hilbert space with scalar product 
(conjugate linear w.r.t. the first variable)
$\langle\cdot,\cdot\rangle$, then $C_\H\in\wtilde \B(\H,\H')$ denotes the isomorphism
defined by $(C_\H\, y)(x):=\langle y,x\rangle$. The Hilbert
adjoint of the densely defined linear operator $A$ is then given by $A^*=
C_\H^{-1}\cdot A'\cdot C_\H$. 
\item $\F$ and $*$ denote Fourier transform and convolution
respectively.
\item $\SC$ denotes the space of distributions and $\E$ is the subspace
of distributions with compact support.
\item $H^s(\RE^n)$, $s\in\RE$, is the usual scale of
Sobolev-Hilbert spaces, i.e. $H^s(\RE^n)$ is the space of tempered
distributions with a Fourier transform which is square integrable
w.r.t. the measure with density $(1+|x|^2)^s$. As usual the strong
dual of $H^s(\RE^n)$ will be represented by $H^{-s}(\RE^n)$.
\item $c$ denotes a generic strictly positive constant
which can change from line to line.
\end{itemize}
\section{Singular Perturbations and Boundary Conditions}
Given the injective self-adjoint operator $A:D(A)\subseteq\H\to\H$, we 
denote by $\D$ the Banach space given by the completion of
$D(A)$ with respect to the norm
$$\|\phi\|_{(A)}:=\|A\phi\|_\H\,.$$
As usual $D(A)$ will be treated as a (dense) subset of $\D$ by means of
the canonical embedding ${\mathcal I}:D(A)\to\D$ which associates 
to $\phi$ the set of all
Cauchy sequences converging to $\phi$. \par 
As in the introduction we consider then a continuous 
linear map 
$$\tau:D(A)\to\X\,,$$
$\X$ is a Banach space, and we will suppose that it satisfies,
besides (h1) and (h2), the further hypothesis
$$
\|\tau\phi\|_\X\le c\,\|A\phi\|_\H\,.\eqno(\h3)
$$
By (h3) $\tau$ admits an extension belonging to
$\B(\D,\X)$; analogously $A$ admits an extension belonging to
$\B(\D,\H)$. By abuse
of notation we will use the same simbols $\tau$ and $A$ to denote these
extensions.\par
Let us now take a sequence $\{\epsilon_n\}_1^\infty\subset\RE$ converging to zero. 
By functional calculus one has
$$
\|(-A\cdot
R(i\epsilon_n)-I)\phi\|^2_\H=\int_{\sigma(A)}d\mu_\phi(\lambda)\,\frac{\epsilon_n^2}{\lambda^2+
\epsilon_n^2}
$$
with $\mu_\phi(\left\{0\right\})=0$ since $0\notin\sigma_{\text{\rm }pp}(A)$. Thus 
$$
1\ge\,\frac{\epsilon_n^2}{\lambda^2+
\epsilon_n^2}\,\longrightarrow 0\,,\qquad\text{\rm $\mu_\phi$-a.e.}
$$
and, by dominated convergence theorem,
$$
\H\-\lim_{n\uparrow\infty}\ -A\cdot R( i\epsilon_n)\phi=\phi\,.
$$
So
$\left\{R(i\epsilon_n)\phi\right\}_1^\infty$ is a Cauchy sequence
in $D(A)$ with respect to the norm $\|\cdot\|_{(A)}$. We can therefore
define $R\in \B(\H,\D)$ by
$$
R\phi:=\widehat D\-\lim_{n\uparrow\infty}\,R( i\epsilon_n)\phi\,,
$$
and then $K(z)\in \wtilde \B(\X',\D)$ by
$$
K(z):=zR\cdot G(z)\,.
$$
Alternatively, using (2), $K(z)$ can be defined by  
$$
K(z)\phi:=\widehat D\-\lim_{n\uparrow\infty}\,
\left(G(i\epsilon_n)-G(z)\right)\phi\,.
$$ 
This immediately implies, using (3),
$$
\forall\,w,z\in\rho(A),\qquad \text{\rm
Range}(K(w)-K(z))\subseteq 
D(A)
$$
and
$$
\forall\,w,z\in\rho(A),\qquad K(w)-K(z)=G(z)-G(w)\,.\eqno(4) 
$$
Also note that
$$
-A\cdot K(z)=zG(z)\,.\eqno(5)
$$
\begin{lemma} The map $$\Gamma:\rho(A)\to\wtilde\B(\X',\X)\,,\qquad
\Gamma(z):=\tau \cdot K(z)
$$
satisfies the relations 
$$
\Gamma(z)-\Gamma(w)=(z-w)\,\breve G(w)\cdot G(z)\eqno(6)
$$
and 
$$
J_\X\cdot \Gamma(z^*)=\Gamma(z)'
\,,\eqno(7)
$$
\end{lemma}
\begin{proof}
Since $K(z)$ is the strong limit of $G(\pm i\epsilon_n)-G(z)$, one has 
$$
\forall \,\ell\in\X'\,,\qquad
\Gamma(z)\ell=\lim_{n\uparrow\infty}\,\hat\Gamma_n(z)\ell\,,
$$
where
$$
\hat\Gamma_n(z):\X'\to\X\,,\qquad
\hat\Gamma_n(z):=\tau
\cdot\left(\,\frac{G(i\epsilon_n)+G(-i\epsilon_n)}{2}-G(z)\,\right)\,.
$$
Thus $\Gamma(z)$ satisfies (6) and 
$$\forall\,\ell_1,\ell_2\in\X'\,,\qquad
\ell_1(\Gamma(z^*)\ell_2)=\left(\ell_2(\Gamma(z)\ell_1)\right)^*$$
(which is equivalent to (7)) since $\hat\Gamma_n(z)$ does 
(see [{\bf 19}, lemma 2.2]).
\end{proof} 
Before stating the next lemma we introduce the following
definition:\par\smallskip\noindent
Given
$\phi\in\H$ and $\psi\in\D$, the writing $\phi=\psi$ will mean that $\phi$ is in
$D(A)$ and ${\mathcal I}\phi=\psi$.\smallskip\noindent
\begin{lemma} Given $\phi\in\H$, $z\in\rho(A)$, suppose there exist
$\psi\in\D$ and $Q\in\X'$ such that 
$$
\phi-G(z)Q=\psi+K(z)Q\,.\eqno(8)
$$
Then the couple 
$(\psi, Q)$ is unique and $z$-independent.
\end{lemma}
\begin{proof}
Let $(\psi_1,Q_1)$, $(\psi_2,Q_2)$ both satisfy
(8). Then
$$
G(z)\,(Q_2-Q_1)=(\psi_1-\psi_2)+K(z)\,(Q_1-Q_2)\,.
$$
By (h2) and the definition of $G(z)$ one has Range$\,G(z)\cap D(A)=\left\{0\right\}$ and so $Q_1-Q_2\in$ Kernel$\,G(z)$. But (h1) implies the
injectivity of $G(z)$ (see [{\bf 19}, remark 2.1]). 
Therefore $(\psi_1,Q_1)=(\psi_2,Q_2)$. The
proof is then concluded observing that $z$-independence follows by
(4).
\end{proof}
We now can extend theorem 1 to the case in which $A$ is injective:
\begin{theorem} Let $A:D(A)\subseteq\H\to\H$ be self-adjoint and
injective, let $\tau:D(A)\to\X$ satisfy {\rm (h1)-(h3)}. 
Given $\Theta\in\wtilde \L(\X',\X)$ self-adjoint, let $D(A^\tau_\Theta)$ be 
the set of $\phi\in\H$ for which there
exist $\phireg\in\D$, $Q_\phi\in D(\Theta)$ such that
$$
\phi-G(z)Q_\phi=\phireg+K(z)Q_\phi
$$
and
$$
\tau\phireg=\Theta\, Q_\phi\,.
$$
Then
$$
A^\tau_\Theta:D(A^\tau_\Theta)\subseteq\H\to\H\,,\qquad
A^\tau_\Theta\,\phi:=A\,\phireg\,.
$$
is a self-adjoint operator which
coincides with $A$ on the kernel of $\tau$ and its resolvent is given by
$$
R^\tau_\Theta(z):=R(z)+G(z)\cdot(\,\Theta+\Gamma(z)\,)^{-1}\cdot\breve
G(z)\,,
\qquad z\in 
W_\Theta^-\cup W_\Theta^+\cup\C\backslash\RE\,.
$$
where $$\Gamma(z):=\tau\cdot K(z)\,.$$
\end{theorem}
\begin{proof} For brevity we define 
$$
\phi_z:=\phi-G(z)Q_\phi\,,\qquad \phi\in D(A^\tau_\Theta) 
$$
and
$$
\Gamma_\Theta(z):=\Theta+\Gamma(z)\,.
$$
Then one has
$$
\Gamma_\Theta(z)
Q_\phi=\tau\phireg+\tau\cdot K(z)Q_\phi=\tau\phi_z\,.
$$
Since $\Gamma(z)$ is a bounded operator satisfying (6) and (7), and
$\Theta$ is self-adjoint, by (h1) and [{\bf 19}, prop. 2.1, remark
2.12], $\Gamma_\Theta(z)$ has a
bounded inverse for any $z\in W^-_\Theta\cup
W^+_\Theta\cup\C\backslash\RE$. 
Therefore $$Q_\phi=\Gamma_\Theta(z)^{-1}\cdot\tau\phi_z$$ 
and 
$$
D(A^\tau_\Theta)\subseteq 
\left\{\phi\in\H\, :\, \phi=
\phi_z+G(z)\cdot\Gamma_\Theta(z)^{-1}\cdot\tau\,\phi_z,\,\phi_z\in
D(A)\right\}\,.
$$
Let us now prove the reverse inclusion.\par
Given $\phi\in\H$, 
$$
\phi=\phi_z+G(z)\cdot\Gamma_\Theta(z)^{-1}\cdot\tau\,\phi_z\,,\qquad \phi_z\in
D(A)\,,
$$
we define $\phireg\in\D$ by
$$
\phireg:=\phi_z-K(z)\cdot\Gamma_\Theta(z)^{-1}\cdot\tau\,\phi_z\,.
$$
Thus one has
\begin{align*}
&\tau\phireg=\tau\phi_z-
\tau\cdot K(z)\cdot\Gamma_\Theta(z)^{-1}\cdot\tau\,\phi_z\\
=&\tau\phi_z-
\tau\phi_z+\Theta\cdot\Gamma_\Theta(z)^{-1}\cdot\tau\,\phi_z
=\Theta\, Q_\phi\,,
\end{align*}
with
$$
Q_\phi:=\Gamma_\Theta(z)^{-1}\cdot\tau\,\phi_z\,.
$$
In conclusion
$$
D(A^\tau_\Theta)=
\left\{\phi\in\H\, :\, \phi=
\phi_z+G(z)\cdot\Gamma_\Theta(z)^{-1}\cdot\tau\,\phi_z,\,\phi_z\in
D(A)\right\}\,.
$$ 
Since, by (5),
$$
A\phireg=A\phi_z-A\cdot K(z)Q_\phi=A\phi_z+zG(z)Q_\phi
$$ 
we have 
$$
(-A^\tau_\Theta+z)\phi=(-A+z)\phi_z\,.
$$ 
Thus $A^\tau_\Theta$ coincides with the operator
constructed in [{\bf 19}, thm. 2.1]; therefore, by (h2), this operator is
self-adjoint, has resolvent given by $R^\tau_\Theta(z)$
and is equal to $A$ on the kernel of $\tau$. 
\end{proof}
\section{Singular perturbations of convolution operators}
Let $\Psi\in\E$. By Paley-Wiener theorem we know
that $\F \Psi$ is a smooth function which is, together with its
derivatives of any order, polynomially bounded. 
Then we define the continuous convolution operator 
$$
\Psi*:\SC\to\SC\,,\qquad \phi\mapsto \Psi*\phi\,, 
$$  
and, supposing $\F \Psi$ real-valued, the restriction of this operator to
the dense subspace
$$
D(\tPsi):=\left\{\phi\in\LD\ :\ \Psi*\phi\in\LD\right\}
$$
provide us with the self-adjoint convolution operator
$$
\tPsi:D(\tPsi)\subseteq\LD\to\LD\,,\qquad \tPsi\phi:=\Psi*\phi\,.
$$ 
Evidently $\tPsi$ is injective if and only if the set of real zeroes of
$\F \Psi$ is a null set. From now on we will therefore suppose that $\F
\Psi$ is real and has a null set of real zeroes. This implies, since $(\text{\rm
Range}\,\tPsi)^\perp=
\text{\rm Kernel}\,\tPsi^*$, that $\tPsi$ has a dense range and, using
(h3), the following
lemma becomes then obvious:
\begin{lemma} The map $\tau:D(\tPsi)\to\X$ can be extended to  
$$\DB:=\left\{\vp\in\SC\ :\ \Psi*\vp\in\LD\right\}\,,
$$
by defining
$$
\tau\vp:=\lim_{n\uparrow \infty}\,\tau\phi_n\,,
$$
where $\left\{\phi_n\right\}_1^\infty\subset D(\tPsi)$ is any sequence such that
$$
L^2\-\lim_{n\uparrow\infty}\,\Psi*\phi_n=\Psi*\vp\,.
$$ 
\end{lemma}
Now we will moreover suppose that $\F \Psi$
is slowly decreasing, i.e. (see \cite{[E]}) we will suppose that there exist
$k>0$ such that for any $\zeta\in\RE^n$ we can find a point
$\xi\in\RE^n$ such that
$$
|\zeta-\xi|\le k\log(1+|\zeta|)\,,
$$
$$
|\F \Psi(\xi)|\ge (k+|\xi|)^{-k}\,.
$$
By [{\bf 11}, thm. 1] we know that $\Psi*:\SC\to\SC$ is surjective if and only if
$\F \Psi$
is slowly decreasing. This is certainly true when $\tPsi$ is a
differential operator, i.e. when $\F \Psi$ is a polynomial. The
hypotheses we made on $\F \Psi$ permit us to state the following
\begin{lemma} Given $\Psi$ as above, 
one has the identification 
$$\DC\,\simeq\,\DB/ \sim\,,
$$
where 
$$
\vp_1\sim\vp_2\quad\iff\quad \Psi*\vp_1=\Psi*\vp_2\,.$$
This identification is given by the isometric maps which to the
equivalence class of Cauchy sequences
$[\left\{\phi_n\right\}_1^\infty]\in \DC$ associates the equivalence class of
distributions $[\vp]\in\DB/ \sim$ such that 
$$
L^2\-\lim_{n\uparrow\infty}\,\Psi*\phi_n=\Psi*\vp\,.
$$
\end{lemma}
\begin{proof} 
Given $[\left\{\phi_n\right\}_1^\infty]\in \DC$, the sequence 
$\left\{\Psi*\phi_n\right\}_1^\infty$ is a Cauchy one in $\LD$ and so it
converges to some $f\in\LD$. Then, by [{\bf 11}, thm. 1], there exists 
$\vp\in\DB$ such that
$\Psi*\vp=f$. Conversely let $\vp\in\DB$; since $\tPsi$ has a dense
range there exists a (unique in $\DC$) sequence
$\left\{\phi_n\right\}_1^\infty
\subset D(\tPsi)$ such that $\Psi*\phi_n$ converges in $\LD$ to $\Psi*\vp$.
\end{proof}
Defining 
$$
\SCT:=\SC/\sim
$$
and then the sum of $\phi\in\LD\subset\SC$ plus
$\psi=[\vp]\in\DC\,\simeq\,\DB/ \sim\,\subseteq\SCT$ by
$$
\phi+\psi:=[\phi+\vp]\in\SCT\,,
$$ 
we can introduce the linear operator
$$
G:\X'\to\SCT\,,\qquad G:=G(z)+K(z)\,.
$$
According to lemma
4, theorem 5 and the definition of $G$, for any $\phi\in 
D(\tPsi^\tau_\Theta)$ we can give the
unique decomposition
$$
\phi=\phireg+GQ_\phi\,.
$$
Thus we can define $D(\tPsi^\tau_\Theta)$ as the set of $\phi\in\LD$ for
which there exists $Q_\phi\in D(\Theta)$ such that
$$
\phi-GQ_\phi=:\phireg\in\DC
$$
and 
$$
\tau\phireg=\Theta\, Q_\phi\,.
$$
\begin{lemma} The definition of $G$ is $z$-independent and 
$$
\forall\,\ell\in\X'\,,\qquad G\ell=[G_*\ell]\,,
$$
where
$$
G_*:\X'\to \SC\,,
$$
is any conjugate linear operator such that
$$
-\Psi*G_*\ell=\tau^*\ell\,.\eqno(9)
$$
Here $\tau^*:\X'\to \SC$ is defined by
$$
\tau^*\ell(\varphi):=(\,\ell(\tau\varphi^*)\,)^*\,, \qquad\varphi\in C_0^\infty(\RE^d)\,.$$
\end{lemma}
\begin{proof} $z$-independence is an immediate consequence of (4). 
By the definition of
$G$ and by (5) there follows
$$
-\Psi*G_*\ell=(-{\Psi*}\,+z)G(z)\ell
$$
and the proof is concluded by the relation  
$$
(-{\Psi*}\,+z)\cdot G(z)=\tau^*\eqno(10)
$$
which can be obtained proceeding as in [{\bf 19}, remark 2.4].
\end{proof}
\begin{remark}
By [{\bf 11}, thm. 1], as $\F \Psi$ is slowly decreasing, the equation (9)
is always resoluble; in particular, denoting the fundamental
solution of $-\Psi*$ by $\G$, when the convolution
$\G*\tau^*\ell$ is well defined (e.g. when $\tau^*\ell\in\E$), one has
$$
G_*:\X'\to\SC\,, \qquad G_*\ell=\G*\tau^*\ell\,.
$$
Analogously, denoting the fundamental solution of
$-{\Psi*}\,+z$ by $\G_z$, one has
$$
G(z):\X'\to\LD\,, \qquad G(z)\ell=\G_z*\tau^*\ell\,.
$$
\end{remark}
\begin{remark} Note that $\tPsi\phireg=\Psi*\vp$ and $\tau\phireg=\tau\vp$ for any
$\vp\in\DB$ such that
$\phireg=[\vp]$. Here we implicitly used the extension given in
lemma 6 and the identification
given in lemma 7. This also implies that $\Gamma(z)$ in
lemma 3 can be re-written as
$$\Gamma(z)=\tau\cdot\left(G_*-G(z)\right)\,.$$
By remark 9, when the convolution is well defined, one can also write
$$\Gamma(z)\ell=\tau\left(\left(\G-\G_z\right)*\tau^*\ell\right)\,.$$
\end{remark}
In conclusion, by making use of the previous lemmata and remarks, we can
restate theorem 5 in the following way:
\begin{theorem} Let $\Psi\in\E$ with $\F \Psi$ real-valued, slowly decreasing and
having a null set of real zeroes, let $\tPsi:D(\tPsi)\subseteq\LD\to\LD$,
$\tPsi\phi:=
\Psi*\phi$, let $\tau:D(\tPsi)\to\X$ satisfy
{\rm (h1)-(h3)}. Given $\Theta\in\wtilde \L(\X',\X)$ self-adjoint, 
let $D(\tPsi^\tau_\Theta)$ be the set of $\phi\in\LD$ for
which there exists $Q_\phi\in D(\Theta)$ such that
$$
\phi-G_*Q_\phi=:\vpreg\in\DB
$$
and 
$$
\tau\vpreg=\Theta\, Q_\phi\,.
$$
Then
$$
\tPsi^\tau_\Theta:D(\tPsi^\tau_\Theta)\subseteq\LD\to\LD\,,\qquad
\tPsi^\tau_\Theta\,\phi:=
\Psi*\vpreg\,,
$$
is a self-adjoint operator which
coincides with $\tPsi$ on the kernel of $\tau$ and its resolvent is given by
$$
R^\tau_\Theta(z):=R(z)+G(z)\cdot(\,\Theta+\Gamma(z)\,)^{-1}\cdot\breve
G(z)\,,
\qquad z\in W_\Theta^-\cup W_\Theta^+\cup\C\backslash\RE\,,
$$
where
$$
\Gamma(z):=\tau\cdot(G_*-G(z))\,.
$$
\end{theorem}

\begin{remark} 
The boundary conditions and the operators $\Gamma(z)$ and 
$\tPsi^\tau_\Theta$ appearing in the previous theorem are
independent of the choice of the representative (see lemma 8) $G_*$ entering in the
definition of $\vpreg$. Indeed any different choice will not
change the equivalence class to which $\vpreg$ belongs, and both $\tau$ and
$\tPsi$ do not depend on the representative in such a class (see remark 10).
\end{remark}
\begin{remark} Proceeding as in [{\bf 19}, remark 2.4] one can give the following
alternative definition of $\tPsi_\Theta^\tau$ where only $Q_\phi\in D(\Theta)$ appears:
$$
\tPsi^\tau_\Theta\,\phi:= \Psi*\phi+\tau^*Q_\phi\,.
$$ 
This is an immediate consequence of identity (10).
\end{remark}
\begin{example}
Let us consider the case $A=\Delta:H^2(\RE^n)\to\LD$. Obviously $A$ is an 
injective convolution operator, thus we can apply to it the
previous theorem. 
\par
A Borel set $F\subset\RE^n$ is called a $d$-set, $d\in(0,n]$, if 
$$
\exists\, c_1,\,c_2>0\ :\ \forall\, x\in F,\ \forall\,r\in(0,1),\qquad
c_1r^d\le\mu_d(B_r(x)\cap F)\le c_2r^d\ ,
$$
where $\mu_d$ is the $d$-dimensional Hausdorff measure and $B_r(x)$ is
the closed $n$-dimensional ball of radius
$r$ centered at 
the point $x$ (see
[{\bf 14}, \S 1.1, chap. VIII]). 
Examples of $d$-sets are $d$-dimensional
Lipschitz submanifolds  and 
(when $d$ is not an integer) self-similar fractals of Hausdorff
dimension $d$ (see [{\bf 14}, chap. II, example 2]). Moreover a finite union of
$d$-sets which intersect on a set of zero $d$-dimensional Hausdorff
measure is a $d$-set.   \par
In the case $0<n-d<4$ we take as the linear operator $\tau$ 
the unique continuous surjective (thus (h1) holds true) map  
$$
\tau_F:H^2(\RE^n)\to B^{2,2}_\alpha(F)\,,\qquad\alpha=2-\,\frac{n-d}{2}
$$
such that, for $\mu_d$-a.e. $x\in F$,
$$
\tau_F\phi(x)\equiv \left\{\phi_F^{(j)}(x)\right\}_{|j|< \alpha}
=\left\{\lim_{r\downarrow 0}\,\frac{1}{\lambda_n(r)}
\int_{B_r(x)}dy\,D^j\phi(y)\,\right\}_{|j|<
\alpha}
\,,\eqno(11)
$$
where $j\in{\mathbb Z}^n_+$, $|j|:=j_1+\dots +j_n$,
$D^j:=\partial_{j_1}\cdots\partial_
{j_n}$ and $\lambda_n(r)$ denotes the
$n$-dimensional Lebesgue measure
of $B_r(x)$. We refer to [{\bf 14}, thms. 1 and 3,
chap. VII] for the existence of the map $\tau_F$; obviously it coincides with the usual
evaluation along $F$ when restricted to smooth functions. 
The definition of the Besov-like (actually Hilbert, see remark 1)
space $B^{2,2}_\alpha(F)$ is quite
involved and we will 
not reproduce it here (see [{\bf 14}, \S 2.1, chap. V]). In the case
$0<\alpha<1$ (i.e. $2<n-d<4$)
things simplify and  $B^{2,2}_\alpha(F)$ can
be defined (see [{\bf 14}, \S 1.1, chap. V]) 
as the Hilbert space of $f\in L^2(F;\mu_F)$ having finite norm
$$\|f\|^2_{B^{2,2}_\alpha(F)}:=
\|f\|^2_{L^2(F)}+\int_{|x-y|<1}d\mu_F(x)\,d\mu_F(y)\,\,
\frac{|f(x)-f(y)|^2}{|x-y|^{d+2\alpha}}\ ,
$$
where $\mu_F$ denotes the restriction of the $d$-dimensional
Hausdorff measure $\mu_d$ to the set $F$. When $\alpha>1$ and $F$ is a 
generic $d$-set the functions $\phi_F^{(j)}\in L^2(F;\mu_F)$ are not
uniquely determined by $\phi^{(0)}_F$; contrarily we may then identify 
$\{\phi_F^{(j)}\}_{|j|< \alpha}$ with the single function
$\phi^{(0)}_F$. This is possible when $F$ preserves Markov's
inequality (see [{\bf 14}, \S 2, chap. II]). Sets
with such a property are closed $d$-sets with $d>n-1$ 
(see [{\bf 14}, thm. 3, \S 2.2, chap. II]), a concrete example
being e.g. the boundary of von Koch's snowflake domain in $\RE^2$ (a
$d$-set with $d=\log 4/\log 3$, see \cite{[W]}). If $F$ has some
additional differential
structure then $B^{2,2}_\alpha(F)\simeq H^{\alpha}(F)$, where $H^{\alpha}(F)$
denotes the usual (fractional) 
Sobolev-Slobodecki\u\i\ space. Some known cases where $B^{2,2}_\alpha(F)\simeq H^{\alpha}(F)$ (for any value of $\alpha>0$) are the following:\p
- $F$ is the graph of a Lipschitz function $f:\RE^d\to \RE^{n-d}$ (see
[{\bf 5}, \S 20]);\p
- $F$ is a bounded manifold of class $C_\gamma$,
$\gamma>\min(3,\max(1,\alpha))$, 
i.e. $F$ has an atlas where the transition maps are of 
class $C^k$, $k<\gamma\le k+1$, and have
derivatives of order less or equal to $k$ which satisfy
Lipschitz conditions of order $\gamma-k$ (see \cite{[J]} for
the case $\gamma>\max(1,\alpha)$ and see [{\bf 5}, \S 24] for the case $\gamma>3$ );\p
- $F$ is a connected complete Riemannian manifold with positive injectivity
radius and bounded geometry, in particular a connected Lie group 
(see [{\bf 23}, \S 7.4.5, \S 7.6.1]).\par
Supposing now that $F$ has a compact closure, let $\chi$ be a smooth
function with a compact support $B$ such that 
$\chi=1$ on $F$. Then  by Sobolev's inequality
one has (from now on $n>4$), 
\begin{align*}
&\|\tau_F\phi\|^2_{B^{2,2}_\alpha(F)}
=\|\tau_F\chi\phi\|^2_{B^{2,2}_\alpha(F)}\le c\,\|\chi\phi\|^2_{H^2(\RE^n)}\\
\le &
c\,(\,\|\phi\|^2_{L^2(B)}+\|\nabla\phi\|^2_{L^2(B)}
+\|\Delta\phi\|^2_{L^2(B)}\,)\\
\le &
c\,(\,\|\phi\|^2_{L^{\frac{2n}{n-4}}(B)}
+\|\nabla\phi\|^2_{L^{\frac{2n}{n-2}}(B)}
+\|\Delta\phi\|^2_{L^2(\RE^n)}\,)\\
\le & c\,\|\Delta\phi\|_{L^2(\RE^n)}\,,
\end{align*}
and so $\tau_F$ satisfies
(h3). Moreover, since 
$$
{\mathcal D}(\Delta)=\left\{\vp\in\SC\,:\,\Delta\phi\in\LD\right\}\subseteq
H^2_{\text{\rm loc}}(\RE^n)
$$
and $\bar F$ is supposed to be compact, the extension of $\tau_F$ to 
${\mathcal D}(\Delta)$ is again defined by (11). 
Denoting the dual of $B^{2,2}_\alpha(F)$ by $B^{2,2}_{-\alpha}(F)$
(the space $B^{2,2}_{-\alpha}(F)$ can be explicitly characterized in the
case $0<\alpha<1$ or when $F$ preserves Markov's inequality, see \cite{[JW2]}), 
hypothesis (h2) is equivalent to $\tau'_F\ell\notin L^2(\RE^n)$
for any $\ell\in B^{2,2}_{-\alpha}(F)\backslash\left\{0\right\}$, where
$\tau'_F\ell\in H^{-2}(\RE^n)$ is defined by
$$
\tau'_F\ell(\phi):=\ell(\tau_F\phi)\,.
$$  
Therefore, as the support of $\tau'_F\ell$ is given by $\bar F$, (h2)
is certainly verified when $\bar F$ has zero Lebesgue
measure. Considering the fundamental solution of
$-\Delta$, given by 
$$
\G(x)=\frac{1}{(n-2)\sigma_n}\,\frac{1}{|x|^{n-2}}\,,
$$ 
$\sigma_n$ the measure of the unitary sphere in $\RE^n$, the
convolution $\G*\tau'_F\ell$ is a well defined distribution as $\tau'_F$
is in $\E$. 
Therefore, by lemma 8 we can choose $G_*$ to be the map
$$
G_*:B^{2,2}_{-\alpha}(F)\to\SC\,,\qquad G_*\ell:=\G*\tau^*_F\ell\,.
$$
Thus, by the previous theorem (and remark 13), supposing that the
$d$-set $F$ has a compact closure of zero Lebesgue measure, given any
self-adjoint operator
$\Theta\in\wtilde
\L(B^{2,2}_{-\alpha}(F),B^{2,2}_\alpha(F))$, $\Theta\in\L(B^{2,2}_\alpha(F))$
if one uses the identification $B^{2,2}_{-\alpha}(F)\simeq
B^{2,2}_\alpha(F)$, we have then the self-adjoint
operator $$\Delta^F_\Theta\,\phi:=\Delta\vpreg\equiv \Delta\phi+\tau_F^*Q_\phi\,,$$ 
where 
$$
\phi=\vpreg+\G*\tau^*_FQ_\phi
\,,\qquad\vpreg\in{\mathcal D}(\Delta)\,,\quad 
Q_\phi\in D(\Theta)\,,$$
and
$$
\left\{\lim_{r\downarrow 0}\,\frac{1}{\lambda_n(r)}\int_{B_r(x)}dy\,
D^j\left(\phi
-\G*\tau^*_FQ_\phi\right)(y)\right\}_{|j|< \alpha}=\Theta\,Q_\phi(x)\,,
$$
$|j|=0$ if $F$ preserves Markov's inequality or 
$B^{2,2}_\alpha(F)\simeq H^\alpha(F)$. When $F=M$, 
$M$ a compact Riemannian manifold, a natural choice for $\Theta$ is given by
$\Theta=(-\Delta_{LB})^{-\alpha}$, where $\Delta_{LB}$ denotes
the Laplace-Beltrami operator. The case in which
$A=(\Delta-\lambda):H^2(\RE^3)\to L^2(\RE^3)$, $\lambda>0$ 
(note that here $0\notin\sigma(A)$, thus
theorem 1 directly applies), and $F$ is a plane circle, 
is treated, without giving
boundary conditions, in
\cite{[KKO]} (also see [{\bf 19}, example 3.2] for connections with
Birman-Kre\u\i n-Vishik theory). 
\end{example}

\end{document}